\documentclass[12pt]{amsart}
\usepackage{amsfonts}

\theoremstyle{definition}

\theoremstyle{remark}

\newcommand{\const}{\mathop{\rm const}\limits}

\newcommand{\grad}{\mathop{\rm grad}\limits}

\begin{document}

\begin{center}

{\bf NECESSARY CONDITIONS FOR FRACTIONAL HARDY-SOBOLEV'S INEQUALITIES} \par

\vspace{3mm}
{\bf E. Ostrovsky and L. Sirota}\\

\vspace{2mm}

{\it Department of Mathematics and Statistics, Bar-Ilan University,
59200, Ramat Gan.}\\
e \ - \ mails: galo@list.ru; \ sirota@zahav.net.il \\

\vspace{4mm}

 {\sc Abstract.} \\

\end{center}

{\it In this short article we obtain some necessary conditions for a so-called
fractional Hardy-Sobolev's inequalities in multidimensional case.\par
 We also give some examples to show the sharpness of these inequalities.} \par

\vspace{3mm}

2000 {\it Mathematics Subject Classification.} Primary 37B30,
33K55; Secondary 34A34, 65M20, 42B25.\\

\vspace{3mm}

{\it Key words and phrases:} Lebesgue-Riesz's norm and spaces, rearrangement invariant (r.i.) Grand and ordinary Lebesgue Spaces, exact estimations, operators, dilation method,
domain, weight. \par

\vspace{3mm}

\section{Introduction. Statement of problem. Notations.}

\vspace{3mm}

{\bf A.  Ordinary Hardy-Sobolev's fractional inequalities.}\par
 The following assertion is called Hardy-Sobolev's (ordinary) difference inequality:

 $$
 \left[ \int_D |u(x)|^q \ |x|^{-\mu} \ dx \right]^{1/q} \le K_{HS}(p,q) \times
 $$
 $$
 \left[\int_D \int_D \frac{|u(x)-u(y)|^p \ |x|^{-\alpha(1)} \ |y|^{-\alpha(2)}}
  {|x-y|^{\beta}} \ dx \ dy  \right]^{1/p}. \eqno(1.1a)
 $$
We will write further $ \mu = \lambda d q. $\par
  Here  $ D $ is open domain with positive Lebesgue measure 
  $ \mu(D) = \int_D dx $ in the whole Euclidean space $ R^d; d=1,2,\ldots $
  equipped with ordinary Euclidean norm $ |x|; \ x \in R^d, $  for instance, whole
space $ R^d $ or its half space or unit ball, $ u = u(\cdot) $ is an arbitrary function
  from the class $ C_0^{\infty}(D),
   \ \alpha(1), \alpha(2), \beta, \lambda = \const, \ p,q = \const \in (1,\infty), $
the finite positive (if there exists)
 "constant" $ K_{HS}(p,q) = K_{HS}(p,q;\alpha(1),\alpha(2);d) $ dependent of the
$ p,q;\alpha(1),\alpha(2),\beta,\mu;d $ but not of the function $ u(\cdot). $ \par

 By means of approximation we can assume $ u \in W^{\mu/q,q}, $  especially when we
 investigate the lower bounds for the constants.  \par

 The finiteness of integrals in the left-hand and right-hand sizes in (1.1a) for
 arbitrary function $ u(\cdot) \in C_0^{\infty}(D) $ entrusts the following conditions
 on the constants $ \alpha(1),\alpha(2),\beta,\mu: $

$$
\mu < d, \ \alpha(1) > -d, \ \alpha(2) > -d, \ \beta < 1, \ \alpha(1) + \alpha(2) -
\beta > -d. \eqno(Ca)
$$
 We assume also that

 $$
 \lambda \in (0,1/(2d-1)). \eqno(Cb)
 $$
{\it We will suppose hereafter the conditions (Ca) and (Cb) are satisfied.} \par

 The following generalization of inequality (1.1a)
 is called Hardy-Sobolev's {\it weight} difference inequality:

 $$
 \left[ \int_D |u(x)|^q \ W_{-\mu}(x)
 \ dx \right]^{1/q} \le K_{HS}(p,q; W) \times
 $$
 $$
 \left[\int_D \int_D \frac{|u(x)-u(y)|^p \ W_{\alpha}(x,y) }
  {W_{\beta}(|x-y|)} \ dx \ dy  \right]^{1/p}. \eqno(1.1b)
 $$

 More general case appears if we write instead $ u(x) $ the function $ u(x) - u(0), $
 for instance

 $$
 \left[ \int_D |u(x)-u(0)|^q \ |x|^{-\mu} \ dx \right]^{1/q} \le K_{HS}^{(0)}(p,q) \times
 $$
 $$
 \left[\int_D \int_D \frac{|u(x)-u(y)|^p \ |x|^{-\alpha(1)} \ |y|^{-\alpha(2)}}
 {|x-y|^{\beta}} \ dx \ dy  \right]^{1/p}. \eqno(1.1c)
 $$

\vspace{3mm}

{\bf B. Mixed Hardy-Sobolev's fractional inequalities.}\par

\vspace{3mm}
 The inequality of a view

 $$
 \left[ \int_D|u(x)|^r \ |x|^{-\mu} \ dx \right]^{1/r} \le K_{M;HS}(p,q,r)\times
 $$
 $$
 \left\{\int_D |y|^{\alpha(2)} dy \
 \left[\int_D  \frac{|u(x)-u(y)|^p \ |x|^{\alpha(1)} \ dx}{|x-y|^{\beta}}  \right]^{q/p}
   \right\}^{1/q} \eqno(1.2a)
 $$
 is said to be (ordinary)  Mixed Hardy-Sobolev's fractional inequality.\par
 The {\it weight} version of  Mixed Hardy-Sobolev's fractional inequality has a view

 $$
 \left[ \int_D|u(x)|^r \ W_{-\mu}(x) \ dx \right]^{1/r} \le K_{MHS;W}(p,q,r)\times
 $$
 $$
 \left\{\int_D  dy \
 \left[\int_D  \frac{|u(x)-u(y)|^p \ W_{\alpha}(x,y)
 \ dx}{W_{\beta}(|x-y|)}  \right]^{q/p}  \right\}^{1/q}. \eqno(1.2b)
 $$

 The functions $  W_{\alpha}(\cdot), W_{\beta}(\cdot), W_{-\mu}(\cdot) $ are
 {\it weight} function, i.e. are measurable positive almost everywhere functions. \par
\vspace{3mm}
{\bf C. Hardy-Sobolev's fractional derivative-difference inequalities.}\par
\vspace{3mm}
 By definition, the inequality of a view

 $$
\left[\int_D \int_D \frac{|u(x)-u(y)|^p \ |x|^{-\alpha(1)} \ |y|^{-\alpha(2)}}
 {|x-y|^{\beta}} \ dx \ dy  \right]^{1/p} \le
 $$
$$
K_{DD;HS}(p,s) \times \int_D \left[|\nabla u(x)|^s \ |x|^{-\mu} \ dx \right]^{1/s} \eqno(1.3c)
$$
is called fractional derivative-difference inequality.\par
Here for the vector $ x = \{ x_1,x_2, \ldots, x_d \}  $
$$
\nabla u(x) = \grad u \stackrel{def}{=} ( \partial u/\partial x_1, \partial u/\partial x_2, \ldots, \partial u/\partial x_d).
$$

 We omit here and in the next pilcrow the obvious weight generalization of
these inequalities.\par
\vspace{3mm}
{\bf D.  Surface Hardy-Sobolev's fractional inequalities.}\par
\vspace{3mm}
 The following assertion is named as surface Hardy-Sobolev's difference
 inequality:

 $$
 \left[ \int_S |u(x)|^q \ |x|^{-\mu} \ \sigma(dx) \right]^{1/q} \le K_{S;HS}(p,q) \times
 $$
 $$
 \left[\int_D \int_D \frac{|u(x)-u(y)|^p \ |x|^{-\alpha(1)} \ |y|^{-\alpha(2)}}
  {|x-y|^{\beta}} \ dx \ dy  \right]^{1/p}. \eqno(1.4a)
 $$

 Here $ S $ is smooth (sub)surface of a boundary  $ \partial D $ of a dimensional
 $ m; 1 \le m \le d-1 $ with correspondent surface measure
 $ d \sigma(x) = \sigma(dx). $ \par

\vspace{3mm}
{\bf  Our aim is finding of some necessary conditions for the fractional
Hardy-Sobolev's inequalities.\par
 We obtain also the lower bounds for the constants $ K_{HS}(p,q), K_{M;HS}(p,q,r),\\
 K_{S;HS}(p,q), K_{DD;HS}(p,s) $ and consider some generalizations on the so-called
 Grand Lebesgue spaces instead classical Lebesgue-Riesz's $ L_p $ spaces. }\par
\vspace{3mm}
 Some upper estimations for Hardy-Sobolev's fractional inequalities
 see, e.g. in the works
\cite{Bogdan1}, \cite{Dyda1}, \cite{Dyda2},\cite{Dyda3}, \cite{Frank1},\cite{Frank2},
\cite{Frank3}, \cite{Kufner1},\cite{Kufner2}, \cite{Kufner3}, \cite{Loss1},
\cite{Maligranda1}, \cite{Maligranda2}, \cite{Maz'ya1} etc. \par
 The one-dimensional case $ d=1 $ for the ordinary Hardy-Sobolev's inequality
 was investigated before by  Jakovlev \cite{Jakovlev1}  and Grisvard
 \cite{Grisvard1}.\par
About applications of these inequalities see, e.g. \cite{Frank2},\cite{Frank3},
\cite{Dyda3},\cite{Kufner3}, \cite{Maz'ya1}.\par

\bigskip

 We use the symbols $C(X,Y),$ $C(p,q;\psi),$ etc., to denote positive
constants along with parameters they depend on, or at least
dependence on which is essential in our study. To distinguish
between two different constants depending on the same parameters
we will additionally enumerate them, like $C_1(X,Y)$ and
$C_2(X,Y).$ The relation $ g(\cdot) \asymp h(\cdot), \ p \in (A,B), $
where $ g = g(p), \ h = h(p), \ g,h: (A,B) \to R_+, $
denotes as usually

$$
0< \inf_{p\in (A,B)} h(p)/g(p) \le \sup_{p \in(A,B)}h(p)/g(p)<\infty.
$$
The symbol $ \sim $ will denote usual equivalence in the limit
sense.\par
 We will denote as ordinary the indicator function
$$
I(x \in A) = 1, x \in A, \ I(x \in A) = 0, x \notin A;
$$
here $ A $ is a measurable set.\par
\bigskip

\section{Main result: necessary conditions for the Hardy-Sobolev's fractional inequalities}

\vspace{3mm}

{\bf A.  Ordinary Hardy-Sobolev's fractional inequalities.}\par
\vspace{3mm}

{\bf Theorem 1A.} Let in the inequality (1.1a) $  D = R^d $ and suppose (1.1a) be satisfied
for some non-constant function $ u(\cdot) \in C_0^{\infty}. $ Then

$$
\frac{d-\mu}{q}= \frac{2d+\alpha(1)+\alpha(2)-\beta}{p}. \eqno(2.1a)
$$
{\bf Proof} used the so-called dilation method belonging to  G.Talenty \cite{Talenty1}.
Namely, let $ \theta=\const \in (0,\infty). $ The dilation operator $ T_{\theta} $
may be defined as follows:

$$
u_{\theta}(x) = T_{\theta}u(x) \stackrel{def}{=} u(x/\theta).
$$
 Note that if $ u(\cdot) \in C_0^{\infty}, $ then $ u_{\theta}(\cdot) \in
 C_0^{\infty}.  $ \par
 We obtain substituting the function $ u_{\theta}(\cdot) $ into inequality (1.1a) instead the  function $ u(\cdot) $  after changing variables $ x=\theta z, y = \theta v: $

$$
\theta^{(d-\mu)/q} \ \left[ \int_D |u(x)-u(0)|^q \ |x|^{-\mu} \ dx \right]^{1/q} \le K_{HS}(p,q) \times
$$

$$
\theta^{(2d+\alpha(1)+\alpha(2)-\beta)/p} \
\left[\int_D \int_D \frac{|u(x)-u(y)|^p \ |x|^{-\alpha(1)} \ |y|^{-\alpha(2)}}
 {|x-y|^{\beta}} \ dx \ dy  \right]^{1/p}.
$$

 Since the value $ \theta $ is arbitrary in the set $(0,\infty), $ we conclude
$$
(d-\mu)/q = (2d+\alpha(1)+\alpha(2)-\beta)/p,
$$
Q.E.D.\par
\vspace{3mm}
Analogously may be proved the following results of this section.\par
{\bf B. Mixed Hardy-Sobolev's fractional inequalities.}\par

\vspace{3mm}

{\bf Theorem 1B.} Let in the inequality (1.2a) $  D = R^d $ and suppose (1.2a) be
satisfied for some non-constant function $ u(\cdot) \in C_0^{\infty}. $ Then

$$
\frac{d-\mu}{r}= \frac{d+\alpha(2)-\beta}{p} + \frac{d+\alpha(1)}{q}. \eqno(2.1b)
$$

\vspace{3mm}
{\bf C. Hardy-Sobolev's fractional difference-derivative inequalities }\par

\vspace{3mm}

{\bf Theorem 1C.} Let in the inequality (1.3a) $  D = R^d $ and suppose (1.3a) be satisfied
for some non-constant function $ u(\cdot) \in C_0^{\infty}. $ Then

$$
\frac{d-\mu}{s} - 1 = \frac{2d+\alpha(1)+\alpha(2)-\beta}{p}.\eqno(2.1c)
$$

\vspace{3mm}
{\bf D. Hardy-Sobolev's surface fractional difference-derivative inequalities }\par

\vspace{3mm}

 We consider here the case when

 $$
 x \in S \ \Leftrightarrow x_{m+1} = x_{m+2} =\ldots = x_n = 0.
 $$
We obtain using as before at the same dilation method as when in the proof of
inequality (1.4a) holds then

$$
\frac{m-\mu}{q} = \frac{2d+\alpha(1)+\alpha(2)-\beta}{p}.\eqno(2.1d)
$$

\bigskip

\section{Weight generalizations of Hardy-Sobolev's fractional inequalities.}

\bigskip
 We consider in this section the Hardy-Sobolev's {\it weight} difference inequality
 in the whole space $ D= R^d: $

 $$
 \left[ \int_{R^d} |u(x)|^q \ W_{-\mu}(x)
 \ dx \right]^{1/q} \le K_{HS}(p,q; W) \times
 $$
 $$
 \left[\int_{R^d} \int_{R^d} \frac{|u(x)-u(y)|^p \ W_{\alpha}(x,y) }
  {W_{\beta}(|x-y|)} \ dx \ dy  \right]^{1/p}. \eqno(3.1)
 $$
 A new notations. For some finite constant $ \alpha_0, \beta_0, \mu_0; \
\alpha_{\infty}, \beta_{\infty}, \mu_{\infty} $ we introduce a
functions, to be presumed non-zero and integrable:

$$
\underline
{W}
_{-\mu,0}(z) =
\inf_{\theta \in (0,1)}\frac{W_{\mu}(\theta z)}{\theta^{-\mu_0}},
$$

$$
\underline
{W}
_{\beta,0}(z) =
\inf_{\theta \in (0,1)}\frac{W_{\beta}(\theta z)}{\theta^{\beta_0}},
$$

$$
\overline
{W}
_{\alpha,0}(z,v) =
\sup_{\theta \in (0,1)}\frac{W_{\alpha}(\theta z,\theta v)}{\theta^{\alpha_0}};
$$

\vspace{3mm}

$$
\underline
{W}
_{-\mu,\infty}(z) =
\inf_{\theta \in (1,\infty)}\frac{W_{\mu}(\theta z)}{\theta^{-\mu_{\infty}}},
$$

$$
\underline
{W}
_{\beta,\infty}(z) =
\inf_{\theta \in (1,\infty)}\frac{W_{\beta}(\theta z)}{\theta^{\beta_{\infty}}},
$$

$$
\overline
{W}
_{\alpha,\infty}(z,v) =
\sup_{\theta \in (1,\infty)}\frac{W_{\alpha}(\theta z,\theta v)}{\theta^{\alpha_{\infty}}}.
$$

\vspace{3mm}

{\bf Theorem 3.1.} If the inequality (3.1) is satisfied for any non-zero function
$ u \in C_0^{\infty}(R^d),$ then

$$
\frac{d-\mu_0}{q} \ge \frac{d+\alpha_0 - \beta_0}{p},\eqno(3.2a)
$$

$$
\frac{d-\mu_{\infty}}{q} \le \frac{d+\alpha_0 - \beta_0}{p}, \eqno(3.2b)
$$

 As a {\bf corollary:} when $ \mu_0 = \mu_{\infty} = \mu, \alpha_0 = \alpha_{\infty}=\alpha,
 \beta_0 = \beta_{\infty}=\beta, $ then both the inequalities (3.2a) and (3.2b)
 reduced to the known relation

$$
\frac{d-\mu}{q} = \frac{d+\alpha - \beta}{p}. \eqno(3.2c)
$$

{\bf Proof} is alike to the proof of theorem 2.1; it used the Talenty dilation method and splitting into two cases: $ \theta \in (0,1)  $ and $ \theta \in (1,\infty). $ \par

\bigskip

\section{Lower bounds for constants in Hardy-Sobolev's fractional inequalities.  Examples.}

\vspace{3mm}

We denote by $\overline{K}_{HS}, \overline{K}_{M;HS}, \overline{K}_{S;HS}, \overline{K}_{DD;HS} $
the {\it minimal values } of the constants
$$
K_{HS}= K_{HS}(p,q), K_{M;HS}= K_{M;HS}(p,q,r), K_{S;HS}= K_{S;HS}(p,q),
K_{DD;HS} = K_{DD;HS}(p,s)
$$
in the (correspondingly)
Hardy-Sobolev's, Mixed Hardy-Sobolev's, Surface Hardy-Sobolev's, and
Differential - Difference Hardy-Sobolev's inequalities for the whole space $ D = R^d. $ \par
 Note that for some particular cases of domains $  D $ (half-spaces etc.) the exact values of these constants was calculated by K.Bogdan and B.Dyda \cite{Bogdan1},
R.Frank, R.Seiringer \cite{Frank1}.\par

 For instance, $ \overline{K}_{HS} = \overline{K}_{HS}(p,q) =
 \overline{K}_{HS}(p,q;\alpha(1),\alpha(2),\lambda) = $

$$
\sup_{u \in C_0^{\infty}(R^d)}
 \left\{\left[ \int_D |u(x)|^q \ |x|^{-\mu} \ dx \right]^{1/q}:\\
 \left[\int_D \int_D \frac{|u(x)-u(y)|^p \ |x|^{-\alpha(1)} \ |y|^{-\alpha(2)}}
  {|x-y|^{\beta}} \ dx \ dy  \right]^{1/p} \right\}.\eqno(4.0)
 $$
 Recall that $  \mu = \lambda d p, \lambda \in (0,1). $ \par
{\bf Theorem 4.a}. Suppose the domain $  D  $  contains some  ball

$$
B(t) = \{x: |x| \le t \}, \ B = B(1)
$$
with the center in origin of the whole space $  R^d. $ \par

 Let also in the ordinary Hardy-Sobolev's inequality $ 1 \le p <
d/\lambda $ and $  \alpha(1) = \alpha(2)=0, \  \beta = d(1+\lambda p). $
Assume that the condition (2.1a) is satisfied. Then

$$
C_1(\lambda;d) \left[\frac{p}{|p-d/\lambda|}\right]^{\lambda} \le
\overline{K}_{HS}(p,p;0,0,\lambda) \le
$$

$$
C_2(\lambda;d) \cdot
\left[\frac{p}{|p-d/\lambda|} \right]. \eqno(4.1a)
$$

{\bf Proof.} The upper estimation contains in fact in
\cite{Bogdan1}; see also \cite{Frank1}; see also
\cite{Grisvard1}, \cite{Jakovlev1}, \cite{Kufner2}. \par
Without loss of generality we can assume $ d=1; $ the multidimensional case $ d \ge 2 $
is investigated analogously, with at the same (counter)example.\par
 It is sufficient to consider the one-dimensional case $ d=1$  and
 to obtain the lower estimate to consider the following example
 in the case when the domain $ D $ contains the unit ball of the set $ R^d,  $ with the center in origin:

 $$
 u_0(x) = |\log |x| \ |\cdot I(|x| \le 1).
 $$

 We have consequently  as $ p \in [1,1/\lambda), \ p \to 1/\lambda-0:  $

 $$
 L^p \stackrel{def}{=} 2 \int_0^1 x^{-\lambda p} \ |\log x| \ dx  =
 \frac{\Gamma(p+1)}{|1-\lambda p|^{p+1}};
 $$

$$
L \asymp |p-1/\lambda|^{1+1/p} \asymp |p-1/\lambda|^{1+\lambda}.
$$
 Further,

 $$
 R^p :=  \left[\int_D \int_D \frac{[|\log|x| - \log|y|]^p \ } {|x-y|^{\beta}} \ dx \ dy  \right] \asymp
 \left[\int_B \frac{[|\log|x| - \log|y|]^p \ }{|x-y|^{\beta}} \ dx \ dy  \right] \asymp
$$

$$
\int_0^1 \rho^{d-1 -\beta} \ |\log \rho| \ d \rho \times
\int_{0}^{2 \pi} \frac{|\log |\tan \phi| \ |^p}
{|\cos \phi - \sin \phi|^{\beta} }  d \phi =:I_1\cdot I_2.
$$
 Note that the second integral $ I_2 $ is bounded when $ 1 \le p < 1/\lambda $ and

 $$
 I_1 \sim \frac{C}{|1-\lambda p|}.
 $$

This completes the proof of theorem 4.a.\par

{\bf Remark 4.1a.} We obtain in more general case when $ \alpha(1) \ge 0, \alpha(2) \ge 0, \alpha :=\alpha(1) + \alpha(2) > 0,  \ \mu = \alpha - \lambda d p $ denoting

$$
p_0 = \frac{1}{\lambda} + \frac{\alpha}{\lambda d}:
$$

$$
C_1(\alpha(1),\alpha(2),\lambda;d) \left[\frac{p}{|p-p_0|}\right]^{1/p_0} \le
\overline{K}_{HS}(p,p;\alpha(1),\alpha(2),\lambda) \le
$$

$$
C_2(\alpha(1),\alpha(2),\lambda;d) \cdot
\left[\frac{p}{|p-p_0|} \right].
$$

\vspace{3mm}

 Analogously may be obtained the following results.\par

{\bf Theorem 4.b}. Let in the mixed ordinary Hardy-Sobolev's inequality $ 1 \le r <
d/\lambda $ and $ \beta = d(1+\lambda p). $ Assume the condition (2.1b) is satisfied.
Then

$$
C_3(\alpha(1),\alpha(2), \lambda;d) \left[\frac{1}{|r-d/\lambda|}\right]^{\lambda} \le \overline{K}_{M;HS}(p,q) \le
$$

$$
C_4(\alpha(1),\alpha(2), \lambda;d) \cdot
\left[\frac{1}{|r-d/\lambda|} \right]. \eqno(4.1b)
$$

\vspace{3mm}

{\bf Theorem 4.c}. Let in the differential-difference Hardy-Sobolev's inequality
$ \beta = d(1+\lambda p) $ and let the condition (2.1c) be satisfied. Then

$$
C_5(\alpha(1),\alpha(2), \lambda;d) \le \overline{K}_{DD;HS}(p,s) \le
C_6(\alpha(1),\alpha(2), \lambda;d).  \eqno(4.1c)
$$

\vspace{3mm}

{\bf Theorem 4.d}. Let in the surface Hardy-Sobolev's inequality $ 1 \le q <
m/\lambda $ and $ \beta = d(1+\lambda p). $  Assume also the condition (2.1d) is satisfied. Then

$$
C_7(\alpha(1),\alpha(2), \lambda;d,m) \left[\frac{q}{|q-m/\lambda|}\right]^{\lambda} \le \overline{K}_{S;HS}(p,q) \le
$$

$$
C_8(\alpha(1),\alpha(2), \lambda;d,m) \cdot
\left[\frac{q}{|q-m/\lambda|} \right]. \eqno(4.1d)
$$

\vspace{3mm}

\section{ Generalization on a Bilateral Grand Lebesgue Spaces}

\vspace{3mm}

  We recall briefly the definition and needed properties of these spaces.
  More details see in the works \cite{Fiorenza1}, \cite{Fiorenza2}, \cite{Ivaniec1},
   \cite{Ivaniec2}, \cite{Ostrovsky1}, \cite{Ostrovsky2}, \cite{Kozatchenko1},
  \cite{Jawerth1}, \cite{Karadzov1} etc.

\vspace{3mm}

For $a$ and $b$ constants, $1 \le a < b \le \infty,$ let $\psi =
\psi(p),$ $p \in (a,b),$ be a continuous positive
function such that there exists a limits (finite or not)
$ \psi(a + 0)$ and $\psi(b-0),$  with conditions $ \inf_{p \in (a,b)} > 0 $ and
 $\min\{\psi(a+0), \psi(b-0)\}> 0.$  We will denote the set of all these functions
 as $ \Psi(a,b). $ \par
 The Bilateral Grand Lebesgue Space (in notation BGLS) $  G(\psi; a,b) =
 G(\psi) $ is the space of all measurable
functions $ \ f: D \to R \ $  or $ \ f: R^d \to R \ $  endowed with the norm

$$
||f||G(\psi) \stackrel{def}{=}\sup_{p \in (a,b)}
\left[ \frac{ |f|_p}{\psi(p)} \right], \eqno(5.1)
$$
if it is finite.\par

 The  $ G(\psi) $ spaces over some measurable space $ (X, F, \mu) $
with condition $ \mu(X) = 1 $  (probabilistic case)
appeared in an article \cite{Kozatchenko1}.\par
 The BGLS spaces are rearrangement invariant spaces and moreover interpolation spaces
between the spaces $ L_1(R^d) $ and $ L_{\infty}(R^d) $ under real interpolation
method \cite{Carro1}, \cite{Jawerth1}. \par
It was proved also that in this case each $ G(\psi) $ space coincides
with the so - called {\it exponential Orlicz space,} up to norm equivalence. In others
quoted publications were investigated, for instance,
 their associate spaces, fundamental functions
$\phi(G(\psi; a,b);\delta),$ Fourier and singular operators,
conditions for convergence and compactness, reflexivity and
separability, martingales in these spaces, etc.\par
 Let $ g: X \to R $ be some measurable function such that $ \exists (a,b), 1 \le a < b
 \le \infty, $ such that $ \forall p \in (a,b) \ \Rightarrow |g|_p < \infty. $

We can then introduce the non-trivial function $ \psi_g(p) $ as follows:

$$
\psi_g(p) \stackrel{def}{=} |g|_p, \ p \in (a,b).  \eqno(5.2)
$$
This choosing of the function $ \psi_g(\cdot) $ will be called {\it natural
choosing. }\par

{\bf Remark 1.} If we introduce the {\it discontinuous} function

$$
\psi_r(p) = 1, \ p = r; \psi_r(p) = \infty, \ p \ne r, \ p,r \in (a,b) \eqno(5.3)
$$
and define formally  $ C/\infty = 0, \ C = \const \in R^1, $ then  the norm
in the space $ G(\psi_r) $ coincides with the $ L_r $ norm:

$$
||f||G(\psi_r) = |f|_r.
$$

Thus, the Bilateral Grand Lebesgue spaces are the direct generalization of the
classical exponential Orlicz's spaces and Lebesgue spaces $ L_r. $ \par

The BGLS norm estimates, in particular, Orlicz norm estimates for
measurable functions, e.g., for random variables are used in PDE
\cite{Fiorenza1}, \cite{Ivaniec1}, theory of probability in Banach spaces
\cite{Ledoux1}, \cite{Kozatchenko1},
\cite{Ostrovsky1}, in the modern non-parametrical statistics, for
example, in the so-called regression problem \cite{Ostrovsky1}.\par

\vspace{3mm}

 Let us introduce the following linear operators acting on the function
 $ u(\cdot)\in C_0^{\infty}(R^d): $

$$
\delta_{\lambda}[u](x,y) = \frac{u(x)-u(y)}{|x-y|^{\lambda d}},\eqno(5.4)
$$

$$
S_{\lambda}[u](x) = \frac{u(x)}{|x|^{\lambda d}}. \eqno(5.5)
$$
 Let also $ \psi_2 = \psi_2(p) $ be some function from the class
 $ \Psi(a,b; R^d \times R^d) $ relative a new {\it potential} measure

 $$
 \nu(dx,dy) = \frac{dx dy}{|x-y|^{\lambda d}}\eqno(5.6)
 $$
 and such that  $ b \le 1/\lambda $ or conversely $ a \ge 1/\lambda. $\par
 Put

$$
\psi_1(p) = \overline{K}_{HS}(p) \cdot \psi_2(p).
$$

\vspace{3mm}

{\bf Theorem 5.1a.}

\vspace{3mm}

$$
||S_{\lambda}[u]||G\psi_2 \le 1 \cdot
||\delta_{\lambda}[u] \ |x-y|^{-\lambda d}||G\psi_1, \eqno(5.7)
$$
where the constant "1" is the best possible.\par

\vspace{3mm}

{\bf Proof} is very simple. Let
$$
\delta_{\lambda}[u] \ |x-y|^{-\lambda d} \in G\psi_2;
$$
without loss of
generality we can suppose $ ||\delta_{\lambda}[u]  |x-y|^{-\lambda d}||G\psi_2= 1.  $\par

 From the direct definition of the norm in Grand Lebesgue spaces it follows

 $$
 \left[\int_{R^d} \int_{R^d}
 \frac{|\delta_{\lambda}[u]|^p(x,y) }{|x-y|^{\lambda d}} dx dy \right]^{1/p} \le
 \psi_2(p), \ p \in (a,b).
  $$
We obtain using the Hardy-Sobolev's inequality  with $ q = p: $

$$
\left[\int_{R^d}|S_{\lambda}[u]|^p(x) \ dx \right]^{1/p} \le  \overline{K}_{HS}(p) \cdot
\psi_2(p)= \psi_1(p),
$$
which is equivalent to the assertion of our theorem.\par
 The precision of the constant " 1 " follows immediately from the main result of
 paper \cite{OSIR}.\par
 Note that it follows from upper estimation for the constant $ K_{HS}(p) $
that if we define a new function $ \psi_3(p)  $ as follows:
$$
\psi_3(p) = \cdot \frac{p \ \psi_2(p)}{|1/\lambda - p|},
$$
then

$$
||S_{\lambda}[u]||G\psi_3 \le C \ \cdot ||\delta_{\lambda}[u]||G\psi_1.
$$
 This result is weakly exact in the following sense. Let $ \psi_4(p) $ be each function
 from the class $ G\Psi(a,b), $ where either $ a = 1/\lambda $ or $ b =1/\lambda $
 for which

 $$
 \lim_{p \to 1/\lambda} \frac{p \ \psi_4(p)}{|p-1/\lambda|^{\lambda}}=0.
 $$
 Then for the function $ u_0(x) = |\log |x| \ | \cdot I(|x|\le 1) $

$$
\lim_{p \to 1/\lambda}
\frac{||S_{\lambda}[u_0]||G\psi_4}
{||\delta_{\lambda}[u_0]||G\psi_1} = \infty.
$$

 \vspace{3mm}
  The {\it mixed, or anisotropic Grand Lebesgue Spaces} was introduced in \cite{OSIR2}. Indeed, let  $ u = u(x), \ x \in R^n $ be measurable function: $ u:R^n \to R.$ Recall
 that the anisotropic Lebesgue space $ L_{\vec{p}} $ consists on all the functions
 $ f $ with finite norm

 $$
 |f|_{\vec{p}}\stackrel{def}{=}
\left( \int_{R^m_1} \mu_1(dx_1) \left( \int_{R^{m_2}}\mu_2(dx_2) \ldots
 \left( \int_{R^{m_l}}
|f(\vec{x})|^{p_1} \mu_l(dx_1) \right)^{p_2/p_1} \right)^{p_3/p_2} \ldots \right)^{1/p_l}. \eqno(5.8)
 $$

  Here $ m_j = \dim x_j, \ \sum_j m_j = n. $ \par
  Note that in general case

 $$
 |f|_{p_1,p_2} \ne |f|_{p_2,p_1},
 $$
but
$$
|f|_{p,p} = |f|_p.
$$
 Observe also that if $ f(x_1,x_2) = g_1(x_1)\cdot g_2(x_2) $ (condition of factorization), then

 $$
 |f|_{p_1,p_2} = |g_1|_{p_1} \cdot |g_2|_{p_2},
 $$
(formula of factorization). \par

 Let $ \nu = \nu(\vec{p}) $ be some continuous positive on the set
 $ Q; \ \vec{p} \in Q $ function such that

$$
  \inf_{p \in Q} \nu(p) > 0, \ \nu(p) = \infty, \ p \notin Q. \eqno(5.9)
$$
We denote the set all of such a functions as $ \Psi_Q. $ \par
 The (multidimensional, anisotropic)  Grand Lebesgue Spaces $ GLS = G_Q(\nu) =G_Q\nu $
space consists on all the measurable functions $ f: R^n \to R $ with finite norms

$$
||f||G_Q(\nu) \stackrel{def}{=} \sup_{\vec{p} \in Q}
\left[ |f|_{\vec{p}} /\nu{\vec{p}} \right]. \eqno(5.10)
$$
 In the further considered case $ n = 2 d, \ \mu_1(dx)=\mu_2(dx) = dx. $\par
{\bf Theorem 5.2.} Let for some $ \nu \in G\nu $

$$
\frac{\delta_{\lambda}(x,y)}{|x-y|^d} \in G\nu(R^d \times R^d).
$$

 Define the domain $ R_r, \ r \ge 1 $ (sub-domain in the plane $ R^2$) by the
 following way:

 $$
 R_r = R_r(\alpha(1),\alpha(2),\beta,\mu;d) =
 \{(p,q): p,q \ge 1, \
  \frac{d-\mu}{r}=\frac{d+\alpha(2)-\beta}{p} +
 \frac{d+\alpha(1)}{q} \} \eqno(5.11)
 $$
and the function

 $$
 \psi_5(r) = \inf_{(p,q) \in R_r} \left[\nu(p,q) K_{M;HS}(p,q) \right].\eqno(5.12)
 $$
Assertion:

 $$
 ||S_{\lambda} u||G\psi_5 \le ||\delta_{\lambda}u \cdot |x-y|^{-\lambda d}||G\nu. \eqno(5.13)
 $$

{\bf Proof.} Without loss of generality we can and do suppose $ ||f||G_Q(\nu)=1; $ then

$$
||\delta_{\lambda} u||_{p,q} \le \nu(p,q).
$$

 It follows from the definition of the norm in BGLS spaces and thee mixed norm inequality (2.1b) that

 $$
 |S_{\lambda} u|_r \le \nu(p,q) K_{M;HS}(p,q), \ (p,q) \in R_r, \eqno(5.14)
 $$
therefore

$$
 |S_{\lambda} u|_r \le
 \inf_{(p,q) \in R_r} \left[\nu(p,q) K_{M;HS}(p,q) \right] = \psi_5(r). \eqno(5.15)
 $$

\end{document}